\documentclass[12pt]{article}

\usepackage{latexsym,amssymb}
\usepackage{graphicx}
\usepackage{graphics}
\usepackage[all]{xy}
\CompileMatrices

\newcommand{\mz}{\mathbb Z}

\newcommand{\mcp}{\mathbb C \mathbb P}
\newcommand{\low}{\mbox{\it low}}

\def\dem{\noindent{\it Proof: \ }}
\def\QED{\hfill$\Box$}

\newtheorem{Teo}{Theorem}[section]
\newtheorem{Lema}[Teo]{Lemma}
\newtheorem{Coro}[Teo]{Corollary}
\newtheorem{Propo}[Teo]{Proposition}
\newtheorem{Remark}[Teo]{Remark}

\newtheorem{Conje}[Teo]{Conjeture}

\begin{document}

\title{\Large\bf
On the $ku$-homology of certain classifying spaces.
\vspace{-2mm}}

\author{\normalsize\sc
Leticia Z\'arate }

\date{}

\maketitle

\begin{abstract}
We calculate the $ku$-homology of $ {\mz}_{p^n} \times {\mz}_{p}$ and ${\mz}_{p^2} \times {\mz}_{p^2}$.
We prove that for these groups the $ku$-homology contains all
the complex bordism information. We construct
a set of generators of the annihilator of the $ku$-toral class.
These elements also generates the annihilator of the $BP$-toral class.
\end{abstract}

{\small


\noindent {{\it Proposed running head}: { On the $\mz_{p^2}
\times \mz_{p^2}$ Conner-Floyd conjecture.}} }

\section{Introduction.}\label{intro}

It is a well known fact that for any prime number $p$, the $BP$-homology of a $p$-local
space $X$ contains all the complex bordism information of $X$. Indeed, we have the Quillen's
splitting theorem:
$ MU_* (X) \otimes_{\mz} \mz_{(p)} = \left( MU_* \otimes_{\mz} \mz_{(p)}  \right ) \otimes_{BP_*} BP_*(X).$

In particular we can take $X$ to be the classifying space of a $p$-group of the form $G = \mz_{p^t} \times \mz_{p^n}$.
We conjecture that for this kind of groups, the local connective complex $K$-theory, here denoted by $ku$, contains
all the complex bordism information.

More precisely, we have the algebraic inclusions:
$BP_* \hookleftarrow BP \langle 1 \rangle_* \hookrightarrow ku_*.$
Here $BP \langle 1 \rangle_*$, $BP_*$ and $ku_*$ are the coefficient rings of the respective homology theories.

In this work we prove that for the groups $G_1 = \mz_{p^t} \times \mz_{p}$ and $G_2 = \mz_{p^2} \times \mz_{p^2}$ there
exists a set of generators of the annihilator of the $ku$-toral class that are elements of the subring $BP\langle 1 \rangle_*$.
G. Nakos \cite{Nak} proved that this elements are also a generating set of the annihilator of the $BP_*$-toral class.

We prove that the annihilator ideal of the $ku_*$-toral class in $ku_* (\mz_{p^2} \wedge \mz_{p^2})$ is generated by:
$$
(p^2, p v^{p-1}, v^{(p-1)(p+2)}).
$$
In \cite{Nak}, G. Nakos proved that the $BP_*$-annihilator of the $BP$-toral class of the group $\mz_{p^2} \times \mz_{p^2}$
is generated by:
$$
(p^2, pv_1, v_1^{p+2}).
$$
In this sense, the $ku$-homology contains all the complex bordism information.

\section{Preliminaries.}\label{s2}

Let $ku$ be the connective complex $K$-theory spectrum and let $v$ be a generator of
$$
\pi_2 (ku) \simeq \mz.
$$

If we consider the canonical orientation:
$$
x = \frac{(1-L)}{v} \in ku^2 (\mcp^{\infty}),
$$
where $L$ is the tautological bundle, the resulting Formal Group Law $F$ is given by:
\begin{equation}\label{FGL}
x +_F y = x + y - vxy.
\end{equation}

For a fixed prime number $p$ we consider the $p$-localization of $ku$ (here also denoted by $ku$) and denote
by $v$ and $x$ the images of the canonical generator and the canonical orientation. We have:
$$
\pi_* (ku) = \mz_{(p)} [v].
$$
In this local spectrum $ku$, the element $x$ determines a Formal Group Law, here also denoted by $F$, which is again given by the formula (\ref{FGL}).

For any natural number $n$ we have the formal power series (that is in fact a polynomial):
$$
[p^n](x) = \sum \limits_{k=0}^{p^n-1} a_{k,n} x^{k+1}
$$
where $a_{k,n} = {p^n \choose k+1} v^k$. We will omit the dependence on $n$ of these coefficients in the
notation since it will be clear from the context.

From here we will use the reduced version of $ku_*$ without further comments. We will denote by
$ku_*(\mz_{p^t} \wedge \mz_{p^n})$ the reduced $ku$-homology of the group $\mz_{p^t} \times \mz_{p^n}$.

It is a well known fact that for any prime number $p$ and any natural number $n$, the group $ku_* (\mz_{p^n})$
is generated by elements
$$
e_i \in ku_{2i - 1} (\mz_{p^n})
$$
for $i \geq 1$ subject to the relations imposed by the $[p^n]$-series. That is,
$$
ku_*(\mz_{p^n}) = \bigoplus_{i \geq 1} ku_* e_i \left / \langle [p^n]\rangle \right .
$$
Here $\langle [p^n] \rangle$ is the $ku$-submodule of $\bigoplus_{i \geq 1} ku_* e_i$ generated by:
$$
\left \{ \sum \limits_{k=0}^{p^n - 1} a_k e_{i-k} \, | \, e_j =0 \,\,{\rm for}\,\, j \leq 0 \right \}.
$$
From the definition of the $a_k$'s we note that those elements with $k+1$ a power of $p$ are of special interest since
$$
\nu_p  {p^n \choose k+1} = n - \nu_p(k+1),
$$
where $\nu_p$ is the usual $p$-valuation, therefore it is natural to have a special notation for these elements.
For any a natural number $i$ we define
$$
g_i = p^i -1;
$$
so $a_{p^i - 1}$ will be denoted by $a_{g_i}$.

\section{Annihilator of the toral class.}\label{s3}

The element $e_{1,n} \in ku_*(\mz_{p^n})$ (the bottom class) is the so called toral class. In the group
$ku_*(\mz_{p^n} \wedge \mz_{p^t})$ we also have a toral class $\tau$, that comes from the
canonical map $\mz^2 \rightarrow \mz_{p^n} \times \mz_{p^t}$.

We have the K\"unneth map
$$
\kappa \colon ku_*(\mz_{p^n}) \otimes_{ku_*} ku_*(\mz_{p^t}) \longrightarrow ku_* (\mz_{p^n} \wedge \mz_{p^t}).
$$
The image of the product of the toral classes $e_{1,n} \otimes e_{1,t}$ under this map is $\tau$. Since this map is injective
we have:
$$
ann_{ku_*} (e_{1,n} \otimes e_{1,t}) = ann_{ku_*} (\tau).
$$

One of the goals of this work is to calculate $ann_{ku_*}(\tau)$ for the groups
$\mz_{p^t} \times \mz_{p}$ and $\mz_{p^2} \times \mz_{p^2}$. In \cite{Nak}, G. Nakos
proved that the $BP_*$-annihilator of the $BP$-toral class of the group $\mz_{p^t} \times \mz_{p}$
is generated by:
$$
(p, v_1^n).
$$
Nakos also proved that the $BP_*$-annihilator of the $BP$-toral class of the group $\mz_{p^2} \times \mz_{p^2}$
is generated by:
$$
(p^2, pv_1, v_1^{p+2}).
$$

Recall that the coefficient ring $BP \langle 1 \rangle_* = \mz_{(p)}
[v_1]$ is a subring of $ku_* = \mz_{(p)} [v]$, with $v_1 = v^{p-1}$
(since the spectrum $ku$ is a wedge sum of $p-1$ suspended copies of
$BP\langle 1 \rangle$). When $p=2$ we have an equality $BP  \langle
1 \rangle= ku$ (as homology theories), and we have proved in
\cite{tesis} that the $ku$-annihilator of the toral class of the
group $\mz_{4} \times \mz_{4}$ is given by:
$$
(4, 2v, v^4)
$$
We have by Nakos' work that it is a also a set of generators of the $BP$-annihilator of $BP$-toral class of $\mz_4 \times \mz_4$.

On the other hand $ku_*(\mz_{p^t} \wedge \mz_{p^n})$ has a direct sum decomposition, since
we have the Landweber split short exact sequence:
\begin{equation}\label{Land}
0 \rightarrow ku_*(\mz_{p^t}) \otimes_{{}_{ku_*}} ku_*(\mz_{p^n})  \rightarrow ku_*(\mz_{p^t} \wedge \mz_{p^n}) \rightarrow
\sum {}{\rm Tor}_1^{ku_*} (ku_*(\mz_{p^t}),ku_*(\mz_{p^n})) \rightarrow 0.
\end{equation}

Since as we will see in the next section we can approximate the first and the third term of this short exact sequence
by an spectral sequence; the annihilator of the toral class contains all the information of the group
$ku_*(\mz_{p^t} \wedge \mz_{p^n})$. This is the reason because we say that $ku$ contains all the complex
bordism information when the $ku$-annihilator of the $ku$-toral class has a set of generators that also generates
the $BP$-annihilator of the $BP$-toral class.

\section{The spectral sequence.}\label{s4}

In this section we construct a spectral sequence that will give,
up to extensions, the group $ku_*(\mz_{p^t} \wedge \mz_{p^n})$. We have a conjecture
about the behavior of this spectral sequence in the general case; in this work
we prove the conjecture for the cases $t>1$, $n=1$ and $t=2=n$.

We fix some notation. Let $F$ the free $ku_*$-module in generators $\alpha_i$'s for $i \geq 1$, that is,
$$
F = \bigoplus_{i \geq 1} ku_* \alpha_i.
$$

For fixed natural numbers $n \geq 1$ and $t \geq n$ consider the map
$$
\begin{array}{rcl}
\partial_t \colon  F \otimes_{ku_*} ku_*(\mz_{p^n})& \longrightarrow & F \otimes_{ku_*} ku_*(\mz_{p^n}) \\
\\
\alpha_i \otimes e_j & \longrightarrow & \sum \limits_{k=0}^{p^t - 1} a_k \alpha_{i-k} \otimes e_j
 \end{array}
$$
where $\alpha_h = 0$ if $h \leq 0$. Here $a_k \in ku_*$ are the coefficients of the $[p^t]$-series.
Note that
$$
{\rm coker} (\partial_t) = ku_*(\mz_{p^t}) \otimes_{ku_*} ku_*(\mz_{p^n}),
$$
since $\partial_t$ imposes the relations of the $[p^t]$-series in the first factor.

We consider the chain complex:
\begin{equation}\label{complex}
\cdots \longrightarrow 0 \longrightarrow F \otimes_{ku_*} ku_*(\mz_{p^n}) \longrightarrow F \otimes_{ku_*} ku_*(\mz_{p^n}) \longrightarrow 0 \longrightarrow \cdots
\end{equation}
where the only non trivial map is given by $\partial_t$.

Note that every element of the module $F \otimes_{ku_*} ku_*(\mz_{p^n})$ has a unique expression modulo $p^n$, therefore we can define a
filtration:
\begin{itemize}
\item $|c| = 0$ for $c \in ku_*$.
\item $|\alpha_i| = p^t +1$.
\item $|e_j| = p^{t-1} + 1$.
\end{itemize}

We will denote by $c [i,j] = c \, \alpha_i \otimes e_j$ for any $c \in ku_*$.

\begin{Conje}\label{conje}
In the spectral sequence associated to the chain complex {\rm (\ref{complex}),} with the defined filtration,
there exist only $n$ families of differentials given by:
$$
p^k [i,j] \longrightarrow p^{n-k-1} v^{h(k)} [i-g_k, j -(t-n+1)g_{k+1} + (t-n)g_k].
$$
Here $h(k) = (t-n+1)g_{k+1} - (t-n-1) g_k$ and $0 \leq k \leq n-1$.
\end{Conje}

In this work we will prove this conjecture for the cases $n=1$, $t>1$ and $n=2=t$.
We have significative advances on the proof of this Conjecture
in the cases $n=2$, $t>2$ and $n=t$ for a natural number $n>2$. The main problem in the
general case is to prove that the elements:
$$
p^k [i,j] \,\,\, {\rm with} \,\,\, i \leq g_k \,\,\, {\rm or} \,\,\, j \leq (t-n+1) g_{k+1} -(t-n)g_k
$$
are permanent cycles. This problem does not appear in the case $n=1$. For the case
$n=2=t$ we prove by hand calculation that the element $[g_2 + g_1 -1, g_1]$ is a
permanent cycle, and as we will see it will be sufficient.

\section{Relations in $F \otimes_{ku_*} ku_*(\mz_{p^2})$.}\label{s5}

In this section we will prove some technical Lemmas that will be the tools to prove the Conjecture~\ref{conje}
in the cases we are dealing with in this work. We begin defining the algebraic Smith morphisms.

For any ordered pair of natural numbers $(i,j)$ with $i,j \geq 0$ we have the Smith morphism $\phi_{i,j}$
given by:
$$
\xymatrix
{
F \otimes_{ku_*} ku_*(\mz_{p^n}) \ar[rr]^{\phi_{i,j}} & & F \otimes_{ku_*} ku_*(\mz_{p^n})\\
[a,b] \ar[rr] & & [a-i,b-j].
}
$$
Note that these morphisms are compatible with $\partial_t$, that is, we have a commutative diagram:
$$
\xymatrix
{
F \otimes_{ku_*} ku_*(\mz_{p^n}) \ar[rr]^{\phi_{i,j} } \ar[dd]_{\partial_t}& & F \otimes_{ku_*} ku_*(\mz_{p^n}) \ar[dd]^{\partial_t} \\
\\
F \otimes_{ku_*} ku_*(\mz_{p^n}) \ar[rr]^{\phi_{i,j} } & & F \otimes_{ku_*} ku_*(\mz_{p^n}).
}
$$
As a consequence we have that for any element
$S= \sum \limits_{i \geq 0} c_i \, [a_i, b_i] \,\, \in \,\, {\rm Im} \,\, \partial_t$
and any ordered pair of natural numbers $(i,j)$ the element $S_{i,j}=\phi_{i,j}(S)$ is
also in ${\rm Im} \,\, \partial_t$. Moreover if
$$
\partial_t \left ( \sum \limits_{k \geq 0} c_k \, [a_k',b_k'] \right ) = S
$$
with higher filtration term given by $c_0 \, [a_0',b_0']$, then
$$
\partial_t \left ( \sum \limits_{k \geq 0} c_k \, [a_k'-i,b_k'-j] \right ) = S_{i,j}.
$$
That is, $S_{i,j}$ is the image under $\partial_t$ of
lower filtration terms than $c_0 \, [a_0',b_0']$.

\begin{Lema}\label{1}
In $F \otimes_{ku_*} ku_*(\mz_{p^n})$, the elements $p^n [a,b]$ with
$1 \leq b \leq g_1$ are zero.
\end{Lema}

\dem We proceed by induction on $b$. Recall that the $ku_*$-module $ku_*(\mz_{p^n})$ is given by:
$$
\bigoplus_{i\geq 1} ku_* e_i \left / \langle [p^n] \rangle \right . ,
$$
therefore for $b=1$ we have that
$$
p^n [a,1] = p^n \alpha_a \otimes e_1 = \alpha_a \otimes (p^n e_1) = 0,
$$
since the second factor is a relation in $ku_*(\mz_{p^n})$ because $a_0 = p^n$. Suppose
we have proved the assertion for $1 \leq b < g_1$, therefore we have
$$
p^n [a, b+1] = - \sum \limits_{i=1}^{b} a_i [a,b+1-i].
$$
Note that the coefficient $a_i$ is divisible by $p^n$ for $0 \leq i < g_1$; this
implies, by the inductive hypothesis, that $p^n[a,b+1]=0$. \QED

For the $[p^2]$-series we have that:
\begin{center}
\begin{tabular}{ll}
$a_i = w_i p^2 v^i$ &for $i+1$ not divisible by $p$. \\
$a_{kp-1} = u_k p v^{kp -1}$ & for $k \leq p$.
\end{tabular}
\end{center}
Here $w_i$ and $u_k$ are units in $ku_*$, with $w_0 = 1$ and $u_p = 1$.

In some results we will need to have total knowledge of the coefficients in the formulas we are dealing with, at least for the ``{\it high}'' filtration terms
and for those terms that are not $p$-divisible. For this
purposes we define the following polynomials with coefficients in $\mz_{(p)}$.

For $0 \leq k \leq p-2$ we define $q_k(x_1, \ldots, x_k) \in \mz_{(p)}[x_1, \ldots, x_k]$ as:
$$
q_0 = - 1 \,\,\,\, {\rm and}  \,\,\,\, q_k(x_1, \ldots, x_k) = -\sum \limits_{i=0}^{k-1} x_{k-i} q_i(x_1, \ldots, x_i) \,\,\,\, {\rm for} \,\,\,\, k \geq 1.
$$
For $0 \leq k \leq p-2$ and $1 \leq n \leq p-2$ we define $q_k^{(n)}(x_1, \ldots, x_k) \in \mz_{(p)}[x_1, \ldots, x_k]$ as:
$$
q_k^{(n)}(x_1, \ldots, x_k) = \left \{
\begin{array}{lcl}
-\sum \limits_{i=0}^{n} x_{k-i} q_i(x_1, \ldots, x_i) & {\rm if} & n \leq k-2 \\
\\
q_k(x_1, \ldots, x_k) & {\rm if} & n \geq k-1.
\end{array} \right .
$$
We will denote by $q_k^{(n)} = q_k^{(n)} (w_1, \ldots, w_k)$.

\begin{Lema}\label{2}
For $2 \leq k \leq p+1$ the following equality holds in $F \otimes_{ku_*} ku_*(\mz_{p^2})$.
$$
p^2 [a,kg_1] = \sum \limits_{i=0}^{p-2} u_1 q_i p v^{g_1+i} [a,(k-1)g_1-i]
+ \sum \limits_{t=0}^{(k-2)g_{{}_1}} z_{t,k} p v^{2g_1 +t} [a,(k-2)g_1 - t].
$$
Here $z_{t,k} \in ku_*$ are coefficients possibly multiples of $p$.
\end{Lema}

\dem We proceed by induction. For $k=2$, we will prove by induction that for $1 \leq n \leq p-3$ the following equality holds:
\begin{equation}\label{rela}
p^2 [a,2g_1] =  - \sum \limits_{i=n+1}^{p-2} q_i^{(n)} v^i p^2 [a,2g_1-i] + \sum \limits_{i=0}^{n} u_1 q_i p v^{g_1+i} [a,g_1-i]
\end{equation}
From the relation imposed in the second factor, we have:
$$
p^2 [a,2g_1] = -\sum \limits_{i=1}^{p-2} w_i p^2 v^{i} [a, 2g_1 - i] - u_1 p v^{g_1} [a,g_1].
$$
The first sumand of the previous equality gives:
$$
-w_1 p^2 v [a,2g_1 - 1] = \sum \limits_{i=1}^{p-3} w_1 w_i p^2 v^{i+1} [a, 2g_1 - (i+1)] + u_1 w_1 p v^{g_1 + 1} [a,g_1 - 1].
$$
Note that the last terms of the relations above don't appear by Lemma~\ref{1}. Mixing both relations we obtain:
$$
\begin{array}{rl}
p^2 [a,2g_1] & = \sum \limits_{i=2}^{p-2} (w_1 w_{i-1} - w_i) p^2 v^i [a, 2g_1 - i] - u_1 p v^{g_1} [a,g_1] + u_1 w_1 p v^{g_1 +1} [a, g_1 -1] \\
\\
             & = - \sum \limits_{i=2}^{p-2} q_i^{(1)} p^2 v^i [a, 2g_1 - i] + \sum \limits_{i=0}^{1} u_1 q_i p v^{g_1 + i} [a,g_1 - i]
\end{array}
$$
Suppose that the equality (\ref{rela}) has been proved for $1 \leq n \leq p-4$, therefore $p^2 [a,2g_1] $ is equal to:
$$
\begin{array}{l}
- \sum \limits_{i=n+2}^{p-2} q_i^{(n)} p^2 v^i [a, 2g_1 - i] + \sum \limits_{i=0}^{n} u_1 q_i p v^{g_1 + i} [a,g_1 - i] \\
\\
+ \sum \limits_{i=1}^{g_{{}_1} -n-2} w_i q_{n+1} p^2 v^{n+i+1} [a, 2g_1 - (n+i+1)]
+ u_1 q_{n+1} v^{g_1+n+1} [a,g_1-(n+1)].
\end{array}
$$
Therefore $p^2 [a,2g_1] =- \sum \limits_{i=n+2}^{p-2} q_i^{(n+1)} p^2 v^i [a, 2g_1 - i] + \sum \limits_{i=0}^{n+1} u_1 q_i p v^{g_1 + i} [a,g_1 - i]$.
So we have that:
$$
p^2[a, 2g_1] = -q_{p-2} p^2 [a, g_1 +1] + \sum \limits_{i=0}^{p-3} u_1 q_i p v^{g_1 + i} [a,g_1 - i]
=  \sum \limits_{i=0}^{p-2} u_1 q_i p v^{g_1 + i} [a,g_1 - i],
$$
and the result follows for $k=2$.

Now we suppose that the result is valid for $2 \leq k \leq p$. From the relation imposed in the second factor
we obtain:
\begin{equation}\label{r1}
\begin{array}{rl}
p^2[a,(k+1)g_1] = & -\sum \limits_{i=1}^{p-2} w_i p^2 v^i [a, (k+1)g_1 -i] - u_1 p v^{g_1} [a,kg_1] \\
\\
& -\sum \limits_{i=p}^{2p-2} w_i p^2 v^i [a, (k+1)g_1 - i] - u_2 p v^{2p-1} [a,(k-1)g_1 - 1] \\
\\
&\,\,\,\,\,\,\, + \cdots +  \\
\\
& -\sum \limits_{i=(k-1)p}^{kp-2} w_i p^2 v^i [a, (k+1)g_1 - i] - u_k p v^{kp-1} [a,p-k].
\end{array}
\end{equation}
Note that when $k=p$ the last summand does not appear since the second coordinate is zero; also note that
when $k<p$ the terms corresponding to $i > kp -1$ are zero by Lemma~\ref{1}.
We have that the first element in the previous equality is given by $-w_1 p^2v[a,(k+1)g_1-1]$,
and it is equal to:
\begin{equation}\label{r2}
\begin{array}{rl}
 &\sum \limits_{i=1}^{p-2} w_1w_i p^2 v^{i+1} [a, (k+1)g_1 -(i+1)] + u_1w_1 p v^{p} [a,kg_1-1] \\
\\
+& \sum \limits_{i=p}^{2p-2} w_1w_i p^2 v^{i+1} [a, (k+1)g_1 - (i+1)] + u_2w_1 p v^{2p} [a,(k-1)g_1 - 2]  + \cdots + \\
\\
+ &\sum \limits_{i=(k-1)p}^{kp-2} w_1w_i p^2 v^{i+1} [a, (k+1)g_1 - (i+1)] + u_kw_1 p v^{kp} [a,p-(k+1)].
\end{array}
\end{equation}
We apply the inductive formula to the element $w_1 w_{p-2} p^2 v^{g_1}[a,kg_1]$ that appears in (\ref{r2}).
Also we apply this formula (up to a Smith morphism) to the elements $-w_i p^2 v^i [a,(k+1)g_1 - i]$ that appears in (\ref{r1})
and to the corresponding elements $w_1w_i p^2 v^{i+1} [a,(k+1)g_1 - (i+1)]$ that appears in (\ref{r2}). Here $p \leq i \leq 2p -2$.
Since all the elements in this formulas are divisible by $p$, replacing (\ref{r1}) in (\ref{r2}) we obtain:
$$
\begin{array}{rl}
p^2[a,(k+1)g_1] =& \sum \limits_{i=2}^{p-2} (w_{i-1} w_1 - w_i) p^2 v^{i} [a, (k+1)g_1 -i] - u_1 p v^{g_1} [a,kg_1] \\
\\
  + & u_1 w_1 p v^p [a,kg_1 - 1] + \sum \limits_{t=0}^{(k-1)g_{{}_1}} y_t p v^{2g_1 + t} [a, (k-1)g_1 - t].
\end{array}
$$
Using the same inductive argument that we applied for the case $k=2$, we can prove that for $1 \leq n \leq p-3$
$$
\begin{array}{rl}
p^2 [a,(k+1)g_1] =&  - \sum \limits_{i=n+1}^{p-2} q_i^{(n)} v^i p^2 [a,(k+1)g_1-i] + \sum \limits_{i=0}^{n} u_1 q_i p v^{g_1+i} [a,kg_1-i] \\
\\
+ &\sum \limits_{t=0}^{(k-1)g_{{}_1} } y_t p v^{2g_1 + t} [a, (k-1)g_1 - t].
\end{array}
$$
Here the coefficients $y_t$ are not necessarily the same in different formulas.
Note that we also applied the inductive formula for $p^2 [a,k g_1]$ (up to a Smith morphism).
Therefore we have:
$$
\begin{array}{rl}
p^2 [a,(k+1)g_1] =& - q_{p-2} p^2 v^{p-2} [a, kg_1 +1] + \sum \limits_{i=0}^{p-3} u_1 q_i p v^{g_1+i} [a,kg_1-i] \\
\\
& + \sum \limits_{t=0}^{(k-1)g_{{}_1}} y_t p v^{2g_1 + t} [a, (k-1)g_1 - t].
\end{array}
$$
Finally we obtain that the element $p^2[a,(k+1)g_1]$ is equal to:
$$
\begin{array}{l}
\sum \limits_{i=0}^{p-2} u_1 q_i p v^{g_1+i} [a,kg_1-i] + \sum \limits_{t=0}^{(k-1)g_{{}_1}} z_{t,k+1} p v^{2g_1 + t} [a, (k-1)g_1 - t].
\end{array}
$$
and the Lemma is proved. \QED

Note that in the previous Lemma we did not take care of the coefficients $z_{t,k+1}$, since as we will see, the ``{\it low}'' filtration and the $p$-divisibility of the terms involved in the second sum will be enough to manipulate them.

\begin{Lema}\label{3}
In the $ku_*$-module $F \otimes_{ku_*} ku_*(\mz_{p^2})$ the following equality holds.
$$
\begin{array}{rl}
p^2[a, g_1 + g_2]  = & - u_1 p v^{g_1} [a,g_2] + \sum \limits_{i=1}^{p-2} u_1 q_i p v^{g_1 + i} [a, g_2 -i] \\
\\
& + \sum \limits_{t=0}^{pg_1} z_t p v^{2g_1+t} [a, pg_1 -t] + \sum \limits_{k=0}^{p-2} q_k v^{g_2 + k} [a, g_1 - k].
\end{array}
$$
Here $q_i = q_i(w_1, \ldots, w_i)$ and $z_t \in ku_*$ is a coefficient possibly multiple of $p$.
\end{Lema}

\dem It is not difficult to verify that $p^2[a,g_1+g_2]$ is equal to:
$$
\begin{array}{rl}
& \sum \limits_{k=2}^{p-2} (-w_k + w_1 w_{k-1}) p^2 v^k [a, g_1+g_2 -k] - u_1 p v^{g_1} [a,g_2] + w_1 u_1 p v^{p} [a, g_2-1] \\
\\
 + & w_1 w_{p-2} p^2 v^{g_1} [a, g_2] - w_p p^2 v^{p}[a, g_2-1] + \sum \limits_{k=p+1}^{2p-2} (-w_k + w_1 w_{k-1}) p^2 v^k [a, g_1+g_2 -k]\\
\\
 + &(-u_2 + w_1 w_{2g_1} p) p v^{2p-1} [a, pg_1 - 1] + (-w_{2p} p + w_1 u_{2} ) p v^{2p} [a, pg_1 - 2] + \cdots +  \\
\\
 + &\sum \limits_{k= pg_{{}_1} +1}^{g_{{}_2}-1}  (-w_k + w_1 w_{k-1}) p^2 v^{k} [a, g_2+g_1-k] - v^{g_2}[a,g_1] + w_1 v^{p^2} [a,g_1-1]\\
\\
 + & w_1 w_{g_2 -1} p^2 v^{g_2} [a,g_1].
\end{array}
$$
Note that the last term is zero by Lemma~\ref{1}. Now we apply Lemma~\ref{2} to the element
$w_1 w_{p-2} p^2 v^{g_1}[a,g_2]$ and (up to Smith morphisms) to the expression:
$$
- w_p p^2 v^{p}[a, g_2-1] + \sum \limits_{k=p+1}^{2p-2} (-w_k + w_1 w_{k-1}) p^2 v^k [a, g_1+g_2 -k]
$$
and we obtain that $p^2[a,g_1 + g_2]$ is equal to:
\begin{equation}\label{ind}
\begin{array}{rl}
&\sum \limits_{k=2}^{p-2} (-w_k + w_1 w_{k-1}) p^2 v^k [a, g_1 + g_2 -k] - u_1 p v^{g_1} [a,g_2] + u_1 w_1 p v^p [c,g_2 -1] \\
\\
+ &\sum \limits_{t=0}^{pg_{{}_1}} z_t p v^{2g_1 + t} [a,pg_1 -t] - v^{g_2} [a, g_1] + w_1 v^{p^2} [a,g_1-1].
\end{array}
\end{equation}

We will prove by induction that for $1 \leq n \leq p-2$ the following equality holds:
$$
\begin{array}{rl}
p^2[a, g_1+ g_2] =&- \sum \limits_{k=n+1}^{p-2} q_k^{(n)} p^2v^k [a,g_1+g_2-k] + \sum \limits_{i=0}^{n} u_1 q_i p v^{g_1 + i} [a,g_2-i] \\
\\
&+ \sum \limits_{t=0}^{pg_{{}_1}} z_t p v^{2g_1 + t} [a, pg_1 -t] + \sum \limits_{k=0}^{n} q_k p v^{g_2 + k} [a,g_1-k].
\end{array}
$$

Since $q_0 = -1$ and $q_1 = w_1$, the beginning of the inductive process is just the formula (\ref{ind}), because
$q_k^{(1)} = w_k - q_1 w_{k-1} = w_k - w_1 w_{k-1}$.

Suppose the assertion has been proved for $2 \leq n \leq p-3$. We have that the element
$q_{n+1}^{(n)} p^2 v^{n+1}  [a, g_1+g_2-(n+1)]$ is equal to:
$$
\begin{array}{l}
 - \sum \limits_{i=1}^{p-2} q_{n+1} w_i p^2 v^{n+i+1} [a, g_1+g_2-(n+i+1)]
- u_1 q_{n+1} p v^{n+p} [a,g_2 - (n+1)] \\
\\
 - \sum \limits_{i=p}^{2p-2} q_{n+1} w_i p^2 v^{n+i+1} [a, g_1+g_2-(n+i+1)]
-u_2 q_{n+1} p v^{2p+n} [a,pg_1 - (n+2)] \\
\\
+ \,\,\,\, A - q_{n+1} v^{p^2+n} [a, g_1 - (n+1)].
\end{array}
$$
Here $A = \sum \limits_{i=2p}^{p^2-2} c_i p v^{n+i+1} [a,g_1+g_2-(n+i+1)]$, where $c_i \in ku_*$ are coefficients possibly multiple
of $p$. Now note that each term in the sum:
$$
\sum \limits_{i=p}^{2p-2} q_{n+1} w_i p^2 v^{n+i+1} [a, g_1+g_2-(n+i+1)]
$$
is a Smith morphism image of $p^2[a,g_2]$ and we can apply Lemma~\ref{2}.
On the other hand, we have for $k \geq n+1$ the equality:
$- q_k^{(n)} + w_{k-(n+1)} q_{n+1} = - q_k^{(n+1)},$
therefore the assertion is valid. So we have:
$$
p^2[a,g_1+g_2] = \sum \limits_{i=0}^{p-2} u_1 q_i p v^{g_1 + i} [a,g_2 - i]
+\sum \limits_{t=0}^{pg_{{}_1}} z_t p v^{2g_1+t} [a, pg_1 -t]
+\sum \limits_{k=0}^{p-2} q_k v^{g_{{}_2} +k} [a, g_1 - k].
$$ \QED

\begin{Remark}\rm
Note that in the previous Lemma we take care of the coefficients of the first and the third sums since the terms of the first one has ``{\it high}'' filtration and
we don't know if the terms of the third one are $p$-divisible.
\end{Remark}

\section{Calculating the differentials.}\label{s6}

In this section we will consider the map
$\partial \colon F \otimes_{ku_*} ku_*(\mz_{p^2}) \longrightarrow F \otimes_{ku_*} ku_*(\mz_{p^2})$
given by the $[p^2]$-series:
$$
\partial ([i,j]) = \sum \limits_{k=0}^{p^2 -1} a_k [i-k,j].
$$
Here the $a_k \in ku_*$ are the coefficients of the $[p^2]$-series and we will use the filtration defined in section~\ref{s4}.
Since each element of $F \otimes_{ku_*} ku_*(\mz_{p^2})$ has a unique expression modulo $p^2$, we define:
$$
|c \, [i,j]| = i(p^2+1) + j (p+1)
$$
for any non zero element $c \in ku_*$.

\begin{Propo}\label{pro}
We have up to units:
\begin{itemize}
\item[a)] $\partial ([1,p]) = p v^{g_1} [1,1].$
\item[b)] $\partial (p[p,p^2] + \,\, {\it low} \,\, ) = v^{g_1+g_2} [1,1].$
\end{itemize}
Here {\it low} stands for lower filtration terms.
\end{Propo}

\dem We have by definition $\partial([1,p]) = p^2 [1,p]$, but
$$
p^2[1,p]  = -\sum \limits_{i=1}^{p-2} w_i p^2 v^i [1,p-i] -u_1 p v^{g_1} [1,1]  = -u_1 p v^{g_1} [1,1].
$$
Since for $1 \leq i \leq p-2$ the term $p^2 v^i [1,p-i]=0$ by Lemma~\ref{1}. This proves a).

We are going to prove that we can find an expression of the form $p[p,p^2] + \,\, {\it low}$ in such a way that its image
under $\partial$ is of the form $A - u_1 v^{g_1 + g_2}[1,1]$. Here $A$ will be a sum of terms $a[i,j]$ divisible
 by $p v^{g_1}$ in such a way
that adding a suitable element on the image of $\partial$ we can replace each $a[i,j]$ with lower filtration terms that are also
divisible by $p v^{g_1}$. We will be careful with the filtration of the terms whose images are used to replace the elements
$a[i,j]$.

We assert that $p v^{g_1} [p,j] \in \,{\rm Im}\, \partial$ if $j \leq g_1^2$.
Indeed, since $j \leq g_1^2$ then $j+g_1 \leq pg_1$, therefore by Lemma~\ref{2} we have:
$$
\begin{array}{rcl}
\partial([p,j+g_1]) &=&\sum \limits_{k=0}^{p-2} w_k p^2 v^k [p-k,j+g_1] + u_1 p v^{g_1} [1,j+g_1] \\
\\
&=& - u_1 \sum \limits_{k=0}^{p-2} w_k p v^{g_1 + k} [p- k,j] + \sum \limits_{t=1}^{j-1} \sum \limits_{k=0}^{p-2} w_k z_t p v^{g_1+k +t} [p-k,j-t] \\
\\
& &+ u_1 p v^{g_1} [1,j+g_1].
\end{array}
$$
We have that $u_1 p v^{g_1} [1,j+g_1]$ is in the image of $\partial$. Indeed,
$\partial([1,j+2g_1])$ is equal to:
$$
-u_1 p v^{g_1} [1,j+g_1] + \sum \limits_{t=1}^{j+g_{{}_1}} z_t p v^{g_1+t} [1,j+g_1-t],
$$
since $j + 2g_1 \leq g_2$. A Smith morphism argument proves that for suitable coefficients $c_t \in ku_*$ we have:
$$
\partial([1,j+2g_1]) + \sum \limits_{t=1}^{j+g_1-1} \partial(c_t [1,j+2g_1-t] ) = -u_1 p v^{g_1}[1,j+g_1].
$$
Therefore we can replace $p v^{g_1} [p,j]$ with elements of the form $c \,p v^{g_1}[i,k]$ with $i \leq p$ and $k \leq j$ and $c \in ku_*$.
Since these elements are Smith morphism images of $p v^{g_1} [p,j]$, we have that $p v^{g_1} [p,j] \in \,{\rm Im}\, \partial$. Note that all
the elements involved in this process are of lower filtration than $[p,p^2]$.

Applying Lemma~\ref{3} (up to a Smith morphism) we have:
$$
p^3[a,p^2]= p \left [ \sum \limits_{i=0}^{p-2} u_1 q_i p v^{g_1 + i} [a,pg_1 -i+1] + \sum \limits_{t=0}^{g_{{}_1}^{2}} z_t p v^{2g_1 +t} [a,g_1^2-t+1] - v^{g_2}[a,1]   \right ].
$$
By Lemma~\ref{2} applied to the first and second sums we obtain:
$$
\begin{array}{rl}
p^3[a,p^2] &= u_1^2 p v^{2g_1} [a, g_1^2 + 1] + \sum \limits_{t=1}^{g_{{}_1}^2} z_t' p v^{2g_1 + t} [a,g_1^2 + 1 -t] \\
\\
             &+ \sum \limits_{t=0}^{(p-2)g_{{}_1}} z_t'' p v^{3g_1 + t} [a, (p-2)g_1 + 1 - t]  - p v^{g_2}[a,1]\\
\\
&  = u_1^2 p v^{2g_1} [a, g_1^2 + 1] + \sum \limits_{t=1}^{g_{{}_1}^2} y_t p v^{2g_1 + t} [a,g_1^2 + 1 -t]  - p v^{g_2}[a,1].
\end{array}
$$
Therefore we have that $\partial (p[p,p^2])$ is equal to:
\begin{equation}\label{d1}
\begin{array}{l}
u_1^2 p v^{2g_1} [p, g_1^2 + 1] + \sum \limits_{t=1}^{g_{{}_1}^2} y_t p v^{2g_1 + t} [p,g_1^2 + 1 -t] - p v^{g_2}[p,1] \\
\\
+ \sum \limits_{i=1}^{p-2} u_1^2 w_i p v^{2g_1 + i} [p-i,g_1^2+1]
+ \sum \limits_{i=1}^{p-2} \sum \limits_{t=1}^{g_{{}_1}^2} w_i y_t p v^{2g_1+ i + t} [p-i,g_1^2 + 1 -t] \\
\\
- \sum \limits_{i=1}^{p-2} w_i p v^{g_2+i}[p-i,1] + u_1 p^2 v^{g_1} [1,p^2] .
\end{array}
\end{equation}
Also we have the equality:
\begin{equation}\label{ultm}
\begin{array}{rl}
u_1 p^2 v^{g_1} [1,p^2] &= -u_1^2 p v^{2g_1} [1, p^2 - g_1] + \sum \limits_{i=1}^{p-2} u_1^2 q_i p v^{2g_1+i} [1, p^2-(g_1 + i)] \\
\\
& + \sum \limits_{t=0}^{g_{{}_1}^2 + 2} u_1 z_t p v^{3g_1 + t} [1,p^2 - (2g_1 + t)] - u_1 v^{g_1 + g_2} [1,1].
\end{array}
\end{equation}
On the other hand $\partial ( u_1 v^{g_1} [p,p^2-g_1])$ is by definition:
\begin{equation}\label{d2}
\begin{array}{rl}
& u_1 v^{g_1} p^2 [p,p^2-g_1] + \sum \limits_{i=1}^{p-2} u_1 w_i p^2 v^{g_1+i} [p-i,p^2 -g_1] + u_1^2 v^{2 g_1} p [1, p^2 - g_1] \\
\\
= & - u_1^2 p v^{2g_1} p [p, p^2-2g_1] -\sum \limits_{i=1}^{p-2} u_1^2 w_i p v^{2g_1 + i} p [p -i, p^2-2g_1] \\
\\
& + \sum \limits_{i=0}^{p-2} \sum \limits_{t=1}^{g_{{}_1}^2} z_{t,i} w_i p v^{2g_1 + t+i} [p-i,g_1^2 -t +1] + u_1^2 p v^{2 g_1}  [1, p^2 - g_1] \\.
\end{array}
\end{equation}

We replace (\ref{ultm}) in (\ref{d1}) and add with (\ref{d2}). We use the fact that the elements of the form
$pv^{g_1} [i,j]$ with $i \leq p$ and $j \leq g_1^2$ are image under $\partial$ of elements of filtration lower than
$p(p^2+1)+p^2(p+1)$ to obtain:
$$
\partial(p[p,p^2] + \, {\it \low} \,) = \sum \limits_{i=1}^{p-2} u_1^2 q_i p v^{2g_1+i} [1, p^2-(g_1 + i)]
+ u_1 z_0 p v^{3g_1} [1,g_1^2+1]- u_1 v^{g_1 + g_2} [1,1],
$$
since $p^2- 2g_1 = g_1^2 + 1$. We have proved that $p v^{g_1} [1,pg_1] \in \,{\rm Im}\, \partial$. On the other hand, we have:
\begin{center}
\begin{tabular}{lcr}
$p^2 -(g_1 +i) \leq pg_1 $ &for& $1 \leq i \leq p-2$
\end{tabular}
\end{center}
and $g_1^2 + 1 < pg_1$. Therefore b) is proved. \QED

\begin{Remark}\rm
Note that the previous proposition gives potential differentials in the spectral sequence we are dealing with:
$$
\begin{array}{rlcl}
d_{\mu_1} \colon & [i,j] &\longrightarrow & p v^{g_1} [i,j-g_1] \\
\\
d_{\mu_2} \colon & p[i,j] &\longrightarrow & v^{g_1 + g_2} [i-g_1, j-g_2].
\end{array}
$$
Here $\mu_1 = g_2$ and $\mu_2 = g_1 (p^2+1) + g_2 (p+1)$. It is not difficult to prove that there can not exist differentials longer than $\mu_2$.
Therefore to prove that $d_{\mu_1}$ and $d_{\mu_2}$ are in fact the only families of differentials in our spectral sequence, it is enough to prove that the element
$[g_1 + g_2 -1,g_1]$ is a permanent cycle.
\end{Remark}

Easy consequences of Lemma~\ref{2} and Lemma~\ref{3} are the following results.

\begin{Coro}\label{*}
For a natural number $0 \leq n \leq p+1$ let denote by $\alpha_n=(p+2-n)g_1$. For $1 \leq k \leq p$ we have that $\partial([ \alpha_k -1, (k+1)g_1 ])$ is equal to:
$$
\begin{array}{l}
-\sum \limits_{i =0}^{p-2} u_1 w_i p v^{g_1+i} [\alpha_k -(i+1),kg_1] +
\sum \limits_{i=0}^{p-2} \sum \limits_{t=1}^{kg_1} z_{t,i} p v^{g_1+t+i} [\alpha_k-(i+1), kg_1-t] \\
\\
+u_1 p v^{g_1} [\alpha_{k+1} -1, (k+1)g_1] + \cdots + u_{k'} p v^{k'p -1} [\alpha_{k+k'}-k',(k+1)g_1] \\
\\
- \sum \limits_{i =k'p}^{\alpha_k -2} u_1 w_i p v^{g_1+i} [\alpha_k-(i+1),kg_1] +
\sum \limits_{i=k'p}^{\alpha_k-2} \sum \limits_{t=1}^{kg_1} w_i z_{t,i} p v^{g_1+t+i} [\alpha_k -(i+1), kg_1-t].
\end{array}
$$
Here $k'$ is such that $k' p  \leq \alpha_k - 1 < (k' + 1)p$.
\end{Coro}

\begin{Coro}\label{**}
The element $\partial([g_1 -1 , g_1 + g_2])$ is equal to:
$$
\begin{array}{l}
 - \sum \limits_{i=0}^{p-3} u_1 w_i p v^{g_1+i} [g_1-(i+1),g_2]
+ \sum \limits_{i=0}^{p-3} \sum \limits_{k=1}^{p-2} u_1 w_i q_k p v^{g_1+i+k} [g_1-(i+1),g_2-k]\\
\\
+ \sum \limits_{i=0}^{p-3} \sum \limits_{k=0}^{pg_{{}_1}} z_{k,i} p v^{2g_1 + i +k } [g_1 - (i+1), pg_1 -k]
+ \sum \limits_{i=0}^{p-3} \sum \limits_{k=0}^{p-2} w_i q_k v^{g_2 +i+k} [g_1-(i+1), g_1 - k].
\end{array}
$$
\end{Coro}

The following results are focused to prove that the element $[(p+2)g_1 -1,g_1]$ is a permanent cycle. The proofs are extensive but easy to follow. We give extended
formulas for the sake of readableness. Most of the formulas look not friendly, but as we will see, many terms can be ignored since they are image under $\partial$
of ``{\it low}'' filtration terms.

\begin{Coro}\label{feo}
We have that $\sum \limits_{k=0}^{p+1} \partial ([(\alpha_n -1,(k+1)g_1])$ is equal to:
$$
\sum \limits_{i=2}^{g_{{}_1}}   \sum \limits_{k=0}^{p-i+1} u_i p v^{ip-1}[\alpha_{k+i} - i, (k+1)g_1] +\sum \limits_{j=1}^{p} A_j +  B_j  + A + B +C  + D
+ v^{g_2} [g_1 -1, g_1].
$$
Here $A_1$ is given by:
$$
A_1 = - \sum \limits_{k=1}^{p} \sum \limits_{i=1}^{p-2} u_1w_i  p v^{g_1 + i} [\alpha_k - (i+1), kg_1],
$$
and for $2 \leq n \leq p$ we have:
$$
A_n = - \sum \limits_{k=1}^{p+2-n} \sum \limits_{i=(n-1)p}^{np-2} u_1w_i  p v^{g_1 + i} [\alpha_k - (i+1), kg_1].
$$
The sum $B_1$ is given by:
$$
B_1 = \sum \limits_{k=1}^{p} \sum \limits_{i=0}^{p-2} \sum \limits_{t=1}^{kg_{{}_1} - 1} z_t w_i p v^{g_1+i + t} [\alpha_k -(i+1), kg_1 -t],
$$
and for $2 \leq n \leq p$ we have that $B_n$ is equal to:
$$
B_n = \sum \limits_{k=1}^{p+2-n} \sum \limits_{i= (n-1)p}^{np-2} \sum \limits_{t=1}^{kg_{{}_1} - 1} z_t w_i p v^{g_1+i + t} [\alpha_k -(i+1), kg_1 -t].
$$
The sums $A$, $B$, $C$ and $D$ are given by:
$$
\begin{array}{l}
A = - \sum \limits_{i=1}^{p-3} u_1 w_i p v^{g_1+i} [g_1-(i+1), g_2].  \\
\\
B= \sum \limits_{i=0}^{p-3} \sum \limits_{t=1}^{p-2} u_1 w_i q_t p v^{g_1+i+t} [g_1 - (i+1), g_2 -t]. \\
\\
C= \sum \limits_{i=0}^{p-3} \sum \limits_{t=0}^{p-2} w_i q_t v^{g_2+i+t} [g_1 - (i+1), g_1 -t].\\
\\
D = \sum \limits_{k=0}^{p-3} \sum \limits_{t=0}^{pg_{{}_1} -1} c_t w_k p v^{2g_1 +t+k} [g_1 -(k+1),pg_1 -t].
\end{array}
$$
Here the elements $c_t \in ku_*$. (They are not necessarily the same in different sums).
\end{Coro}

\dem For $1 \leq k \leq p$ the element $\partial ([a_k-1, (k+1)g_1])$ is given by:
$$
\begin{array}{l}
-\sum \limits_{i =0}^{p-2} u_1 w_i p v^{g_1+i} [\alpha_k -(i+1),kg_1] +
\sum \limits_{i=0}^{p-2} \sum \limits_{t=1}^{kg_1} z_{t,i} p v^{g_1+t+i} [\alpha_k-(i+1), kg_1-t] \\
\\
+u_1 p v^{g_1} [\alpha_{k+1} -1, (k+1)g_1] + \cdots + u_{k'} p v^{k'p -1} [\alpha_{k+k'}-k',(k+1)g_1] \\
\\
- \sum \limits_{i =k'p}^{\alpha_k -2} u_1 w_i p v^{g_1+i} [\alpha_k-(i+1),kg_1] +
\sum \limits_{i=k'p}^{\alpha_k-2} \sum \limits_{t=1}^{kg_1} w_i z_{t,i} p v^{g_1+t+i} [\alpha_k -(i+1), kg_1-t].
\end{array}
$$
The first term of this expression is $-u_1 p v^{g_1} [\alpha_{k}-1, (k+1)g_1]$. On the other hand the first
term of the second row is given by $u_1 p v^{g_1} [\alpha_{k+1}g_1-1, kg_1]$.
When we take the sum over $k$ this terms are canceled, since we
have that the first term in the formula of Corollary~\ref{**} is $-u_1 p v^{g_1} [g_1-1,g_2]$,
and by Lemma~\ref{1} the term $\partial([\alpha_0 -1,g_1])$ is equal to:
$$
u_1 p v^{g_1} [g_2-1,g_1] + \sum \limits_{i=2}^{p-1} u_i v^{ip-1} [p^2-(i-1)p-2,g_1] + v^{g_2} [g_1 -1, g_1],
$$
and the Corollary follows. \QED

Now we are going to prove that the elements in the formula of the previous Corollary that are divisible by
$pv^{g_1}$ and have ``{\it low filtration}'' are in ${\rm Im} \, \partial$. We will be careful with the filtration of the
elements we are using to replace those elements. We want to obtain:
\begin{equation}\label{reduc}
\partial([(p+2)g_1 -1,g_1] + \low) =  v^{g_2} [g_1 -1, g_1] + A_1 + B_1 + A + B + C.
\end{equation}
For to do this, we will use the following result.

\begin{Lema}\label{chop}
For $0 \leq k \leq p-1$ we have, up to units:
$$
\partial ([\alpha_{k+2},(k+2)g_1] + \low) =  p v^{g_1} [\alpha_{k+2}, (k+1)g_1].
$$
\end{Lema}

\dem We will proceed by reverse induction in $k$. For $k = p-1$ we prove that $pv^{g_1} [g_1,b] \in \, {\rm Im}\, \partial$ for $1 \leq b \leq pg_1$.
Indeed, for $b=1$ we have
$$
\begin{array}{rl}
\partial([g_1,p]) &= \sum \limits_{i=0}^{p-2} w_i p^2 v^i [g_1 -i,p] \\
\\
& = - \sum \limits_{i=0}^{p-2} w_i u_1 p v^{g_1+i} [g_1 -i,1].
\end{array}
$$
We will verify that, up to units, $\partial([t,p] + \low) = p v^{g_1} [t,1]$ for $1 \leq t \leq g_1-1$. We have proved in Proposition~\ref{pro}
that $\partial([1,p]) = p v^{g_1} [1,1]$. We suppose that the assertion has been proved for $1 \leq t < g_1 -1$, therefore for $t+1$ we have:
$$
\partial([t+1,p])  = - \sum \limits_{i=0}^{t} w_i u_1 p v^{g_1+i} [t -i+1,1],
$$
by the inductive hypothesis the assertion follows. This proves that $\partial ([g_1,p]+ \low) = p v^{g_1}[g_1,1]$.

Now suppose we have proved that $\partial([g_1,b+g_1] + low) = p v^{g_1} [g_1,b]$ for $1 \leq b < pg_1$. Therefore:
$$
\begin{array}{rl}
\partial([g_1,b+g_1+1]) & = \sum \limits_{i=0}^{p-2} \sum \limits_{t=0}^{p-2} w_i u_1 q_{t} p v^{g_1+i+t} [g_1 -i,b+1-t]\\
\\
& + \sum \limits_{i=0}^{p-2} \sum \limits_{t=0}^{b-g_{{}_1}} w_i z_t p v^{2g_1+t} [g_1-i,b -(g_1+t)+1].
\end{array}
$$
By the inductive hypothesis we have that:
$$
\partial([g_1,b+g_1+1]+low) = -\sum \limits_{i=0}^{p-2}  w_i u_1  p v^{g_1+i} [g_1 -i,b+1].
$$
Applying a Smith morphism to this expression we obtain:
$$
\partial([1,b+g_1+1]+low) = -u_1 p v^{g_1} [1,b+1].
$$
Now is clear that an easy inductive argument proves that $\partial([t,b+g_1+1] + low) = pv^{g_1} [t,b+1]$ for $1 \leq t \leq g_1 -1$. This implies that
$$
\partial([\alpha_{p+1},g_2] + low) = p v^{g_1} [\alpha_{p+1}, pg_1].
$$

Now suppose that for $0 < k \leq p-1$ we have proved:
\begin{equation}\label{induc1}
\partial([\alpha_{k+2},(k+2)g_1]+\low) = p v^{g_1} [\alpha_{k+2}, (k+1)g_1].
\end{equation}
For $k-1$ we have:
$$
\begin{array}{rl}
\partial([\alpha_{k+1},(k+1)g_1]) & = \sum \limits_{i=0}^{p-2} w_i p^2 v^i [\alpha_{k+1},(k+1)g_1] \\
\\
& + u_1 p v^{g_1} [\alpha_{k+2},(k+1)g_1] + \sum \limits_{i=p}^{\alpha_{k+1}} c_i p v^{i} [\alpha_{k+1}-i,(k+1)g_1].
\end{array}
$$
Since for $p \leq i \leq \alpha_{k+1}$ we have $\alpha_{k+2} < \alpha_{k+1}-i$, the inductive hypothesis gives:
$$
\begin{array}{rl}
\partial([\alpha_{k+1},(k+1)g_1]+\low) & = \sum \limits_{i=0}^{p-2} w_i p^2 v^i [\alpha_{k+1}-i,(k+1)g_1] \\
\\
                                       & = \sum \limits_{i=0}^{p-2} \sum \limits_{t=0}^{p-2} w_i u_1 q_t p v^{g_1+i+t}  [\alpha_{k+1}-i,kg_1 - t] \\
                                      \\
                                       & + \sum \limits_{i=0}^{p-2} \sum \limits_{t=0}^{(k-1)g_{{}_1}} w_i z_t p v^{2g_1+i+t}  [\alpha_{k+1}-i,(k-1)g_1 - t].

\end{array}
$$
Now, for $1 \leq i \leq p-2$, it can be proved by reverse induction that:
$$
\partial([\alpha_{k+1}-i, (k+1)g_1] + \low) = p v^{g_1} [\alpha_{k+1}-i, kg_1].
$$
For $i = p-2$ we have by (\ref{induc1}):
$$
\begin{array}{rl}
\partial([\alpha_{k+2} + 1, (k+1)g_1] + \low) & = \sum \limits_{i=0}^{p-2} w_i p^2 v^i [\alpha_{k+2}-i+1,(k+1)g_1] \\
\\
& = \sum \limits_{i=0}^{p-2} \sum \limits_{t=0}^{p-2} u_1 w_i q_t p v^{g_1+i+t} [\alpha_{k+2}-i+1,kg_1 -t] \\
\\
& + \sum \limits_{i=0}^{p-2} \sum \limits_{t=0}^{(k-1)g_{{}_1}} w_i z_t p v^{2g_1+i+t} [\alpha_{k+2}-i+1, (k-1)g_1 -t].
\end{array}
$$
Applying the inductive hypothesis we obtain:
$$
\begin{array}{rl}
\partial([\alpha_{k+2} + 1, (k+1)g_1] + \low) & = \sum \limits_{t=0}^{p-2} u_1 q_t p v^{g_1+t} [\alpha_{k+2}+1,kg_1 -t] \\
\\
& + \sum \limits_{t=0}^{(k-1)g_{{}_1}} z_t p v^{2g_1+t} [\alpha_{k+2}+1, (k-1)g_1 -t].
\end{array}
$$
We can prove inductively for $1 \leq t \leq kg_1$ that:
$$
\partial([\alpha_{k+2}+1, t +g_1] + \low) = p v^{g_1} [\alpha_{k+2}+1,t].
$$
We will only prove the case $t=1$ since in this point the inductive process is clear. We have by induction:
$$
\partial([\alpha_{k+2}+1,p] + \low) = - \sum \limits_{i=0}^{p-2} w_i u_1 p v^{g_1+i} [\alpha_{k+2}-i+1,1],
$$
and by (\ref{induc1}):
$$
\partial([\alpha_{k+2}+1,p] + \low) = -u_1 p v^{g_1} [(p-k)g_1-1,1].
$$
Therefore we have that
$$
\partial([\alpha_{k+2} + 1, (k+1)g_1] + \low) = p v^{g_1} [\alpha_{k+2}+1,kg_1].
$$
Following the same ideas used until here, it is easily proved that for $1 \leq i < p-2$:
$$
\partial([\alpha_{k+1}-i, (k+1)g_1] + \low) = p v^{g_1} [\alpha_{k+1}-i, kg_1].
$$
Therefore we have that:
$$
\begin{array}{rl}
\partial([\alpha_{k+1},(k+1)g_1]+\low) & = \sum \limits_{t=0}^{p-2} u_1 q_t p v^{g_1+t}  [\alpha_{k+1},kg_1 - t] \\
                                      \\
                                       & + \sum \limits_{t=0}^{(k-1)g_{{}_1}} z_t p v^{2g_1+t}  [\alpha_{k+1},(k-1)g_1 - t].

\end{array}
$$
Applying a Smith morphism to the previous expression we have:
$$
\partial([\alpha_{k+1},p] + \low) = -u_1 p v^{g_1} [\alpha_{k+1},1],
$$
and now an easy inductive argument proves that:
$$
\partial([\alpha_{k+1},(k+1)g_1]+\low) = -u_1 p v^{g_1} [\alpha_{k+1},kg_1],
$$
and the Lemma follows. \QED

\begin{Remark}\rm
Note that by the previous Lemma the sum:
$$
\sum \limits_{i=2}^{g_{{}_1}}   \sum \limits_{k=0}^{p-i+1} u_i p v^{ip-1}[\alpha_{k+i} - i, (k+1)g_1]
$$
from the Corollary~\ref{feo} can be ignored. Indeed, all the terms involved in this sum are Smith morphism images of elements of the form $pv^{g_1}[\alpha_{k+2},(k+1)g_1]$ for some $0 \leq k \leq p-1$. Similar arguments can be used for the sums $A_n$, $B_n$ for $n > 2$ and for the sum $D$. Note that for the sum $B_1$ the terms corresponding to $t>p-2$
also can be treated with the previous Lemma, therefore, we can replace $B_1$ in (\ref{reduc}) with:
$$
B_1' = \sum \limits_{k=1}^{p} \sum \limits_{i=0}^{p-2} \sum \limits_{t=1}^{g_{{}_1} - 1} z_t w_i p v^{g_1+i + t} [\alpha_k -(i+1), kg_1 -t].
$$
\end{Remark}

Until here we have proved that:
$$
\partial([\alpha_0 - 1, g_1] + \low) = v^{g_2} [g_1 -1, g_1] + A_1 + B_1' + A + B + C.
$$

\begin{Lema}\label{chopper}
Let denote by $\Delta  = v^{g_2} [g_1 -1, g_1] + A + B + C.$ We have that:
$$
\begin{array}{rl}
\Delta & - \sum \limits_{i=1}^{p-3} \partial (q_i v^i [g_1-(i+1),g_2+g_1]) \\
\\
& = \sum \limits_{i=1}^{p-2} u_1 q_i p v^{g_{{}_1} + i} [g_1-1, g_2-i] + \sum \limits_{i=1}^{p-2}
q_i v^{g_{{}_2} + i} [g_1-1, g_1-i] \\
\\
&+\sum \limits_{i=1}^{p-3} \sum \limits_{k=0}^{p-i-3}
\sum \limits_{t=0}^{pg_{{}_1}} c_t p v^{2 g_{{}_1}+ t + i} [g_1-(k+i+1), pg_1-t].
\end{array}
$$
\end{Lema}

\dem We assert that for $1 \leq n \leq p-4$ we have the following:
$$
\begin{array}{rcl}
\Delta  &-& \sum \limits_{i=1}^{n} \partial (q_i v^i [g_1-(i+1),g_2+g_1]) \\
\\
&= & - \sum \limits_{i=n+1}^{p-3} u_1 q_i^{(n)} p v^{g_{{}_1} + i} [g_1-(i+1), g_2 + g_1] + S_1 + S_2 \\
\\
& & + \sum \limits_{t=1}^{p-2} \sum \limits_{i=n+1}^{p-3} u_1 q_t q_i^{(n)}
p v^{g_{{}_1} + t + i} [g_1-(i+1), g_2-t] \\
\\
& & + \sum \limits_{t=0}^{p-2} \sum \limits_{i=n+1}^{p-3} q_t q_i^{(n)} v^{g_{{}_2 } + t + i} [g_1-(i+1), g_1-t] \\
\\
& & + \sum \limits_{i=1}^{n} \sum \limits_{k=0}^{p-i-3} \sum \limits_{t=0}^{pg_{{}_1}} c_t p
v^{2 g_{{}_1} + t + k +i} [g_1-(k+i+1),pg_1-t].
\end{array}
$$
We will proceed by induction. Recall that $q_1 = w_1$, therefore for $n = 1$ we apply Corollary~\ref{**} (up to a
Smith morphism) to obtain that $\Delta - \partial(q_1 v [g_1-2,g_2+g_1])$ is equal to:
\begin{equation}\label{primas}
\begin{array}{l}
\Delta  + u_1 w_1 p v^{p} [g_1-2,g_2] \\
\\
+ \sum \limits_{k=1}^{p-4} u_1 w_1 w_k p v^{g_{{}_1}+i +k+1} [g_1 - (k+2), g_2] \\
\\
- \sum \limits_{k=0}^{p-4}  \sum \limits_{i=1}^{p-2} u_1 w_1 w_k q_i p v^{g_{{}_1}+k+1} [g_1 - (k+2), g_2-i] \\
\\
- \sum \limits_{k=0}^{p-4}  \sum \limits_{i=0}^{p-2} w_1 w_k q_i v^{g_{{}_2}+i+k+1} [g_1 - (k+2), g_1-i] \\
\\
- \sum \limits_{k=0}^{p-4}  \sum \limits_{t=0}^{pg_{{}_1}} c_t p v^{2 g_{{}_1}+t+k+1} [g_1 - (k+2), pg_1-t]
\end{array}
\end{equation}
We will denote by $A'$, $B'$ and $C'$ the second, third and fourth rows of this sum, that is:
$$
A' = \sum \limits_{k=1}^{p-4} u_1 w_1 w_k p v^{g_{{}_1} +k+1} [g_1 - (k+2), g_2],
$$
and so on. Now note that the first summand of $A$ (when $i=1$) is canceled with the term $u_1 w_1 p v^{p} [g_1-2,g_2]$
of (\ref{primas}), therefore we will omit it in the next calculation. Adding $A$ with $A'$ we have:
$$
A + A' = \sum \limits_{k=2}^{p-3} u_1 (w_1w_{k-1} - w_k) p v^{g_{{}_1} + k}[g_1 -(k+1), g_2].
$$
The sum  $S_1 = \sum \limits_{i=1}^{p-2} u_1 q_i p v^{g_{{}_1} + i} [g_1-1, g_2-i]$ is the term of $B$ when $k=0$
since $w_0 = 1$; the term of $B$ when $k=1$ is given by:
$$
\sum \limits_{i=1}^{p-2} u_1 w_1 q_i p v^{g_{{}_1} + i +1} [g_1-2, g_2-i].
$$
On the other hand, the summand of $B'$ when $k=0$ is equal to:
$$
-\sum \limits_{i=1}^{p-2} u_1 w_1 q_i p v^{g_{{}_1} + i +1} [g_1-2, g_2-i],
$$
therefore:
$$
B + B' = S_1 + \sum \limits_{k=2}^{p-3}
\sum \limits_{i=1}^{p-2} u_1 (w_k - w_1 w_{k-1}) q_i p v^{g_{{}_1} + i +k} [g_1-(k+1), g_2-i].
$$
We have that $S_2 = \sum \limits_{i=1}^{p-2} q_i v^{g_{{}_2} + i} [g_1-1, g_1-i]$ is the summand of $C$ when $k=0$
(note that the term corresponding to $i=0$ has been canceled with $v^{g_2} [g_1 -1 , g_1]$).
It is easily seen that the summand of $C$ for $k=1$ is canceled with the summand of $C'$ for $k=0$. Therefore we
have that:
$$
C + C' = S_2 +  \sum \limits_{k=2}^{p-3}
\sum \limits_{i=0}^{p-2} (w_k - w_1 w_{k-1}) q_i v^{g_{{}_2} + i +k} [g_1-(k+1), g_1-i].
$$
But by definition $q_k^{(1)} = w_k - w_1 w_{k-1}$. This proves the assertion for $n=1$.
Suppose it has been proved for $1 \leq n < p-4$.
Therefore $\Delta - \sum \limits_{k=1}^{n+1} \partial (q_k v^k [g_1 - (k+1),g_2+g_1])$ is equal to:
$$
\begin{array}{l}
- \sum \limits_{k=n+1}^{p-3} u_1 q_k^{(n)} p v^{g_{{}_1} + k} [g_1-(k+1),g_2] + S_1 + S_2 \\
\\
 + \sum \limits_{i=1}^{p-2} \sum \limits_{k=n+1}^{p-3} u_1 q_k^{(n)} q_i p v^{g_{{}_1} + i +k} [g_1-(k+1), g_2-i]\\
\\
 + \sum \limits_{i=0}^{p-2} \sum \limits_{k=n+1}^{p-3} q_k^{(n)} q_i v^{g_{{}_2} + i +k} [g_1-(k+1), g_1-i] \\
 \\
 + \sum \limits_{k=1}^{n} \sum \limits_{i=0}^{p-k-3} \sum \limits_{t=0}^{pg_{{}_1}}
c_t p v^{2 g_{{}_1} +t+ i +k} [g_1-(k+i+1), pg_1-t]\\
\\
 + \sum \limits_{i=0}^{p-n-4} u_1 w_i q_{n+1} p v ^{g_{{}_1}+i +n+1} [g_1-(i+n+2),g_2] \\
 \\
  - \sum \limits_{k=0}^{p-n-4} \sum \limits_{i=1}^{p-2} u_1 w_k q_{n+1} q_i p v ^{g_{{}_1}+i+k+n+1}
[g_1-(k+n+2),g_2-i]
\end{array}
$$
$$
\hspace{-.5cm}
\begin{array}{l}
- \sum \limits_{k=0}^{p-n-4} \sum \limits_{t=0}^{pg_{{}_1}} c_t p v ^{2 g_{{}_1}+t+k+n+1} [g_1-(k+n+2),pg_1-t]\\
\\
 - \sum \limits_{k=0}^{p-n-4} \sum \limits_{i=0}^{p-2} w_k q_{n+1} q_i v ^{g_{{}_2}+i+k+n+1} [g_1-(k+n+2),g_1-i].
\end{array}
$$
When we add the first row (without $S_1+S_2$) with the fifth row we obtain:
$$
 \sum \limits_{k=n+2}^{p-3} u_1 ( q_{n+1} w_{k-n-1} - q_k^{(n)} ) p v^{g_{{}_1} + k} [g_1-(k+1),g_2].
$$
The sum of the second and the sixth rows is equal to:
$$
\sum \limits_{k=n+2}^{p-3} \sum \limits_{i=1}^{p-2} u_1 q_i (q_k^{(n)} - q_{n+1} w_{k-n-1} ) p v^{g_{{}_1} +i+ k}
[g_1-(k+1),g_2-i].
$$
Finally, adding the third and the last row we have:
$$
\sum \limits_{k=n+2}^{p-3} \sum \limits_{i=0}^{p-2} q_i (q_k^{(n)} - q_{n+1} w_{k-n-1} ) v^{g_{{}_2} +i+ k}
 [g_1-(k+1),g_2-i].
$$
Since by definition $q_k^{(n+1)} = q_k^{(n)} - q_{n+1} w_{k-n-1}$ the assertion is proved. The Lemma follows from Lemma~\ref{3} since $\partial( q_{p-3} v^{p-3} [1,g_2+g_1]) = q_{p-3} p^2 v^{p-3} [1,g_2+g_1]$.
\QED

\begin{Remark}\rm
Note that by the Lemma~\ref{chop} we can ignore the sum:
$$
\sum \limits_{i=1}^{p-3} \sum \limits_{k=0}^{p-i-3}
\sum \limits_{t=0}^{pg_{{}_1}-1} c_t p v^{2 g_{{}_1}+ t + i} [g_1-(k+i+1), pg_1-t].
$$
in Lemma~\ref{chopper}. Therefore until here we have proved that:
$$
\partial([\alpha_0 - 1, g_1] + \low) = S_1 + S_2 + A_1 + B_1'.
$$
At this point the crux of the matter is to deal with $S_1 + S_2$, since we have to deal with elements that are not
divisible by $p$. If we prove that this sum can be replaced by the image under $partial$ of ``{\it low filtration}''
terms we are done. Indeed, a Smith morphism argument will take care of $A_1 + B_1'$.
\end{Remark}

\begin{Lema}\label{filgen}
For $1 \leq i \leq p-2$ we have that
$- \sum \limits_{k=0}^{p-3} \partial(q_i q_k v^{i+k} [g_1-(k+1), g_2+g_1-i])$ is equal to:
$$
\begin{array}{l}
- u_1 q_i p v^{g_1 +i} [g_1-1,g_2-i] + \sum \limits_{t=1}^{p-2} u_1 q_t q_i p v^{g_{{}_1}+t+i} [g_1-1, g_2-(i+t)] \\
\\
+ \sum \limits_{t=0}^{p-i-2} q_i q_t v^{g_{{}_2}+t+i} [g_1-1, g_1-(i+t)] \\
\\
+ \sum \limits_{h=0}^{p-3} \sum \limits_{k=0}^{p-h-3} \sum \limits_{t=0}^{pg_{{}_1}}
z_t p v^{2 g_{{}_1} +t+k+i+h} [g_1-(k+h+1), p g_1 - (i+t)].
\end{array}
$$
\end{Lema}

\dem We have that $\partial (q_i v^i [g_1-1, g_1+g_2-i]) - \partial (w_1 q_i v^{i+1} [g_1-2,g_1+g_2-i])$ is equal to:
$$
\hspace{-1cm}
\begin{array}{l}
- u_1 q_i p v^{g_1 + i} [g_1-1,g_2-1] - \sum \limits_{k=2}^{p-3} u_1 q_i q_k^{(1)} p v^{g_1+i+k} [g_1-(k+1),g_2-i]
+ S_1' \\
\\
+\sum \limits_{k=2}^{p-3}\sum \limits_{t=1}^{p-2} u_1 q_i q_k^{(1)} q_t p v^{g_1+k+t+i} [g_1-(k+1),g_2-(i+t) ]+ S_2'\\
\\
- \sum \limits_{h=0}^{1} \sum \limits_{k=0}^{p-3} \sum \limits_{t=0}^{pg_{{}_1}}
z_t p v^{2 g_{{}_1} +t+k+i+h} [g_1-(k+h+1),  pg_1 - (i+t)] \\
\\
+ \sum \limits_{k=2}^{p-3} \sum \limits_{t=0}^{p-i-2} q_i q_t q_k^{(1)} v^{g_2+t+k+i} [g_1 -(k+1),g_1-(t+i)].
\end{array}
$$
Here
$$
S_1'= \sum \limits_{t=1}^{p-2} u_1 q_t q_i p v^{g_1+t+i} [g_1-1,g_2-(i+t)]  ,
$$
and
$$
S_2'= \sum \limits_{t=0}^{p-i-2} q_i q_t v^{g_2+t+i} [g_1-1,g_1-(t+i)] .
$$

It is easily seen by induction, as in the proof of Lemma~\ref{chopper}, that for $n<p-3$ we have
that $- \sum \limits_{k=0}^{n} \partial(q_i q_k v^{i+k} [g_1-(k+1), g_2+g_1-i])$ is equal to:
$$
\begin{array}{l}
S_1' + S_2' - u_1 q_i p v^{g_1+i} [g_1-1,g_2-i] \\
\\
- \sum \limits_{k=n+1}^{p-3} u_1 q_i q_k^{(n)} p v^{g_1+i+k} [g_1-(k+1), g_2-i]) \\
\\
- \sum \limits_{h=0}^{n} \sum \limits_{k=0}^{p-3-h} \sum \limits_{t=0}^{pg_{{}_1}}
z_t p v^{2 g_{{}_1} +t+k+i+h} [g_1-(k+h+1),  pg_1 - (i+t)] \\
\\
+ \sum \limits_{k=n+1}^{p-3} \sum \limits_{t=1}^{p-2} u_1 q_i q_t q_k^{(n)} p v^{g_1+t+i+k} [g_1-(k+1),g_2-(t+i)]   \\
\\
+ \sum \limits_{k=n+1}^{p-3} \sum \limits_{t=0}^{p-i-2} q_i q_t q_k^{(n)} v^{g_2+t+k+i} [g_1 -(k+1),g_1-(t+i)].
\end{array}
$$
Applying the Lemma~\ref{3} to:
$$
- \partial (q_i q_{p-3} v^{i+p-3} [1,g_2+g_1-i]) = q_i q_{p-3} p^2 v^{i+p-3} [1,g_2+g_1-i],
$$
the result follows easily. \QED

Note that by the Lemma~\ref{chop} we can ignore the sum
$$
\sum \limits_{h=0}^{p-3} \sum \limits_{k=0}^{p-h-3} \sum \limits_{t=0}^{pg_{{}_1}}
z_t p v^{2 g_{{}_1} +t+k+i+h} [g_1-(k+h+1), p g_1 - (i+t)],
$$
in the previous Lemma, so we have that:
\begin{equation}\label{ecred}
\begin{array}{rl}
\partial ([\alpha_0 -1,g_1] + \low) & = A_1 + B_1' + S_2 \\
\\
&+ \sum \limits_{i=1}^{p-2} \sum \limits_{t=1}^{p-2} u_1 q_t q_i p v^{g_{{}_1}+t+i} [g_1-1, g_2-(i+t)] \\
\\
&+ \sum \limits_{i=1}^{p-2} \sum \limits_{t=0}^{p-i-2} q_i q_t v^{g_{{}_2}+t+i} [g_1-1, g_1-(i+t)] \\
\end{array}
\end{equation}
On the other hand, when we consider the last row of the previous expression, the term corresponding to $t=0$
is given by $-S_2$, since $q_0 = -1$, therefore:
$$
\begin{array}{rl}
\partial ([\alpha_0 -1,g_1] + \low) & = A_1 + B_1' \\
\\
&+ \sum \limits_{i_0=1}^{p-2} \sum \limits_{i_1=1}^{p-2-i_0} u_1 q_{i_0} q_{i_1} p v^{g_{{}_1}+i_0+i_1}
[g_1-1, g_2-(i_0+i_1)] \\
\\
&+ \sum \limits_{i_0=1}^{p-3} \sum \limits_{i_1=1}^{p-i_0-2} q_{i_0} q_{i_1} v^{g_{{}_2}+i_0+i_1}
[g_1-1, g_1-(i_0+i_1)]
\end{array}
$$
Now we want to deal with the second and the third rows
of the last expression. Using a Smith morphism argument the following result
is a consequence of Lemma~\ref{filgen}.

\begin{Lema}\label{filgen1}
For $1 \leq i_0 \leq p-2$ and $1 \leq i_1 \leq p-2-i_0$ we have that:
$$
- \sum \limits_{k=0}^{p-3} \partial(q_{i_0} q_{i_1} q_k v^{i_0+i_1+k} [g_1-(k+1), g_2+g_1-(i_0 + i_1)])
$$
is equal to:
$$
\begin{array}{l}
- u_1 q_{i_0} q_{i_1} p v^{g_1 +i_0 + i_1} [g_1-1,g_2-(i_0 + i_1)] \\
\\
+ \sum \limits_{t=1}^{p-2} u_1 q_t q_{i_0} q_{i_1} p v^{g_{{}_1}+i_0+i_1+t} [g_1-1, g_2-(i_0+i_1+t)] \\
\\
+ \sum \limits_{t=0}^{p-(i_0+i_1)-2} q_t q_{i_0} q_{i_1}  v^{g_{{}_2}+i_0+i_1+t} [g_1-1, g_1-(i_0+i_1+t)] \\
\\
+ \sum \limits_{h=0}^{p-3} \sum \limits_{k=0}^{p-h-3} \sum \limits_{t=0}^{pg_{{}_1}}
z_t p v^{2 g_{{}_1} +t+k+i_0+i_1+h} [g_1-(k+h+1), p g_1 - (i_0+i_1+t)].
\end{array}
$$
\end{Lema}

Adding the elements of the previous Lemma to the formula (\ref{ecred}) we obtain:
$$
\begin{array}{rl}
\partial ([\alpha_0 -1,g_1] + \low) & = A_1 + B_1' \\
\\
&+ \sum \limits_{i_0=1}^{p-2} \sum \limits_{i_1=1}^{p-2-i_0} \sum \limits_{i_2=1}^{p-i_0-i_1-2}
u_1 q_{i_0} q_{i_1} q_{i_2} p v^{g_{{}_1}+i_0+i_1} [g_1-1, g_2-(i_0+i_1)] \\
\\
&+ \sum \limits_{i_0=1}^{p-4} \sum \limits_{i_1=1}^{p-i_0-3} \sum \limits_{i_1=1}^{p-i_0-i_1-2}
q_{i_0} q_{i_1} q_{i_2} v^{g_{{}_2}+i_0+i_1+i_2} [g_1-1, g_1-(i_0+i_1+i_2)].
\end{array}
$$

Now we can use an inductive argument to obtain that $\partial ([\alpha_0 -1,g_1] + \low)$
is equal to:
\begin{equation}\label{red2}
\begin{array}{l}
A_1 + B_1' + \sum \limits_{i_0=1}^{p-2} \sum \limits_{i_1=1}^{p-i_0-2} \cdots \sum \limits_{i_{n+1}=1}^{p - I_n-2}
u_1 q_{{}_{I_{n+1}}}  p v^{g_{{}_1}+I_{n+1}} [g_1-1, g_2-I_{n+1}] \\
\\
+ \sum \limits_{i_0=1}^{p-(n+3)} \sum \limits_{i_1=1}^{p-i_0-(n+2)} \cdots \sum \limits_{i_{n+1}=1}^{p-I_{n}-2}
q_{{}_{I_{n+1}}} v^{g_{{}_2}+I_{n+1}} [g_1-1, g_1-I_{n+1}].
\end{array}
\end{equation}
Here $I_j = \sum \limits_{k=0}^{j} i_k$ and $q_{{}_{I_j}} = \prod \limits_{k=0}^{j} q_{i_k}$.

When $n+1 = p-3$ we have:
$$
\begin{array}{rl}
\partial ([\alpha_0 -1,g_1] + \low) & = A_1 + B_1' \\
\\
&+ \sum \limits_{i_0=1}^{p-2} \sum \limits_{i_1=1}^{p-i_0-2} \cdots \sum \limits_{i_{p-3}=1}^{p - I_{p-4}-2}
u_1 q_{{}_{I_{p-3}}}  p v^{g_{{}_1}+I_{p-3}} [g_1-1, g_2-I_{p-3}] \\
\\
& + q_{{}_{I_{p-3}}} v^{g_{{}_2}+I_{p-3}} [g_1-1, g_1-I_{p-3}].
\end{array}
$$
Since the upper limit of the $j$-th sum of the last row in (\ref{red2}) is given by
$$p - I_{j-1} - (n+3-j)$$ and $I_{j-1} \geq j$, we have $p - I_{j-1} - (n+3-j) \leq 1$
when $n+1 = p-3$.

Now we consider:
$$
- \sum \limits_{k=0}^{p-3} \partial(q_{{}_{I_{p-3}}} q_k v^{I_{p-3}+k} [g_1-(k+1), g_2+g_1-I_{p-3}])
$$
to obtain
$$
\begin{array}{l}
- u_1 q_{{}_{I_{p-3}}} p v^{g_1 +I_{p-3}} [g_1-1,g_2- I_{p-3}] \\
\\
+ \sum \limits_{t=1}^{p-2} u_1 q_t q_{{}_{I_{p-3}}} p v^{g_{{}_1}+I_{p-3}+t} [g_1-1, g_2-I_{p-3}+t)] \\
\\
+ \sum \limits_{t=0}^{p-I_{p-3}-2} q_t q_{{}_{I_{p-3}}}  v^{g_{{}_2}+I_{p-3}+t} [g_1-1, g_1-I_{p-3}-t)] \\
\\
+ \sum \limits_{h=0}^{p-3} \sum \limits_{k=0}^{p-h-3} \sum \limits_{t=0}^{pg_{{}_1}}
z_t p v^{2 g_{{}_1} +t+k+h+I_{p-3}} [g_1-(k+h+1), p g_1 - I_{p-3}-t].
\end{array}
$$
Here $I_{p-3} = \sum \limits_{k=0}^{p-3} i_k$ with $1 \leq i_0 \leq p-2$ and $1 \leq i_k \leq p - I_{k-1} - 2$.
Note that the sum in the third row of the previous expression only has one term:
$$
q_{{}_{I_{p-3}}} q_{0} v^{g_2 + I_{p-3}} [g_1-1,1] = - q_{{}_{I_{p-3}}} v^{g_2 + I_{p-3}} [g_1-1,1].
$$
Therefore we have proved that $\partial ([\alpha_0 -1,g_1] + \low)$ is equal to:
$$
A_1 + B_1' + \sum \limits_{i_0=1}^{p-2} \sum \limits_{i_1=1}^{p-i_0-2} \cdots \sum \limits_{i_{p-2}=1}^{p - I_{p-3}-2}
u_1 q_{{}_{I_{p-2}}}  p v^{g_{{}_1}+I_{p-2}} [g_1-1, g_2-I_{p-2}].
$$
Since $I_{p-2} \geq g_1$ we have that $g_2 - I_{p-2} \leq pg_1$ and by Lemma~\ref{chop} we obtain:
$$
\partial ([\alpha_0 -1,g_1] + \low)  = A_1 + B_1'.
$$
Applying a Smith morphism to the previous formula plus an inductive argument it is easily verified that each term in the sum $A_1 + B_1 '$ is image under $\partial$ of a low filtration term. This proves that $[(p+2)g_1-1,g_1]$ is a permanent cycle.

\section{The $ku$-homology of the groups.}\label{s7}

In this section we apply the calculations of the previous sections to identify a set of generators of the
annihilator ideals of the toral classes in $ku_* (\mz_{p^2} \wedge \mz_{p^2})$
and $ku_* (\mz \wedge \mz_{p^n})$.

\begin{Teo}
In the spectral sequence of section~\ref{s4} with $t=2=n$ there exist only two families of
differentials and they are given by:
$$
\xymatrix
{
[i,j] \ar[rr]^{d_{\mu_1}\,\,\,\,} & &p v^{g_1} [i,j-g_1].
}
$$
$$
\hspace{1cm}
\xymatrix
{
p[i,j] \ar[rr]^{d_{\mu_2}\,\,\,\,\,\,\,\,\,\,}& &v^{g_1+g_2} [i-g_1,j-g_2].
}
$$
Here $\mu_1 = g_2$ and $\mu_2 = g_1(p^2+1) + g_2 (p-1)$.
\end{Teo}

\begin{Teo}
The annihilator ideal of the $ku_*$-toral class in $ku_* (\mz_{p^2} \wedge \mz_{p^2})$ is generated by:
$$
(p^2, p v^{g_1}, v^{g_1+g_2}).
$$
\end{Teo}

\begin{Teo}
The group $ku_* (\mz_{p^2} \wedge \mz_{p^2})$ is given (up to extensions) by $E_0 \oplus E_1$, where:
$$
E_0 = \bigoplus \limits_{\stackrel{i \geq 1}{j \geq 1}} ku_* \left / (p^2,pv^{g_1},v^{g_1+g_2}) \langle
[i,j] \rangle \right .
$$
and
$$
\begin{array}{rl}
E_1 = & \bigoplus \limits_{\stackrel{i \geq 1}{g_1 \geq j \geq 1}} ku_* \left / (p^2) \langle
[i,j] \rangle \right . \\
\\
&\bigoplus \bigoplus \limits_{\stackrel{i \geq g_1}{j \leq g_2}} ku_* \left / (p^2) \langle
p [i,j] \rangle \right. \\
\\
&\bigoplus \bigoplus \limits_{\stackrel{1 \leq i \leq g_1}{g_1 \leq j}} ku_* \left / (p^2) \langle
p [i,j] \rangle . \right.
\end{array}
$$
\end{Teo}

It is not difficult to verify that when we consider the spectral sequence of section~\ref{s4} with
$n=1$ and $t>1$ we obtain:

\begin{Teo}
In the spectral sequence there exists only one family of
differentials and it is given by:
$$
[i,j] \longrightarrow v^{t g_1} [i,j-t g_1].
$$
\end{Teo}

\begin{Teo}
The annihilator of the $ku_*$-toral class in $ku_* (\mz_{p^t} \wedge \mz_{p})$ is generated by:
$$
(p, v^{t g_1}).
$$
\end{Teo}

\begin{Teo}
The group $ku_* (\mz_{p^t} \wedge \mz_{p})$ is given (up to extensions) by $E_0 \oplus E_1$, where
$$
E_0 = \bigoplus \limits_{\stackrel{i \geq 1}{j \geq 1}} ku_* \left / (p,v^{t g_1}) \langle
[i,j] \rangle \right .
$$
and
$$
\begin{array}{rl}
E_1 = & \bigoplus \limits_{\stackrel{i \geq 1}{t g_1 \geq j \geq 1}} ku_* \left / (p) \langle
[i,j] \rangle \right .
\end{array}
$$
\end{Teo}

\begin{Teo}{\rm \cite{Nak}}
The $BP_*$-annihilator ideals of the $BP$-toral classes in  the groups
$BP_* (\mz_{p^2} \wedge \mz_{p^2})$
and $BP_* (\mz_{p} \wedge \mz_{p^t})$ are generated by:
$$
(p^2, p v_1, v_1^{p+2}) \,\,\,\,\,\,\,\,\,\, {\it and} \,\,\,\,\,\,\,\,\,\,
(p,v_1^t),
$$
respectively.
\end{Teo}

\begin{Coro}
For the $p$-groups $\mz_{p^2} \times \mz_{p^2}$ and $\mz_{p} \times \mz_{p^t}$ the
$ku$-homology contains all the complex bordism information.
\end{Coro}

\bigskip\medskip
Leticia Z\'arate

CEFyMAP - Universidad Aut\'onoma de Chiapas.

4a. Oriente Norte No. 1428.  Entre 13a. y 14a. Norte.

Col. Barrio La Pimienta

Tuxtla Guti\'errez, Chiapas. C.P. 29000

M\'exico

e-mail: {\tt leticia@math.cinvestav.mx}


\begin{thebibliography}{99}

\bibitem{Nak}{G.~Nakos: On the Brown-Peterson
homology of certain classifying spaces,
Ph.D. Thesis, The Johns Hopkins University, 1985.}

\bibitem{tesis}{L.~Z\'arate: On the $BP \langle n \rangle_*$-homomology
of $\mz_{2^e} \times \mz_{2^e}$,
Ph.D. Thesis, CINVESTAV - IPN, 2007.}

\bibitem{BP-osa}{L.~Z\'arate: {\it On the $BP$-homomology
of $\mz_{2^e} \times \mz_{2^e}$},
Journal of K-theory: K-theory and its Applications to Algebra, Geometry,
and Topology, Volume 3, Issue 03, (June 2009) 409-435.}

\end{thebibliography}
\end{document}